\def\QQ         {{\bf Q}}

\def\ii         {{\rm i}}

\def\dim        {{\rm dim}}
\def\sin        {{\rm sin}}

\def\sinh       {{\rm sinh}}
\def\deg        {{\rm deg}}

\def\Box        {{\bullet}}

\documentstyle[twoside,12pt]{article}
\newtheorem{prop}{Proposition}[section]
\newtheorem{dfn}[prop]{Definition}   
\newtheorem{theo}[prop]{Theorem}
\newtheorem{conj}[prop]{Conjecture}
 
\newtheorem{coro}[prop]{Corollary} 
\newtheorem{lem}[prop]{Lemma}

\title{ Elliptic genera of singular varieties, orbifold 
elliptic genus and chiral de Rham complex.} 

\author{
Lev A. Borisov
\\
\small Department of Mathematics,  Columbia University, New York, NY
10027\\
\small e-mail: lborisov@math.columbia.edu\\
Anatoly Libgober\footnote{Supported by NSF grant. The authors also 
thank Arthur Greenspoon for careful reading of the manuscript.} \\
\small Department of Mathematics, University of Illinois, Chicago, IL
60607\\
\small e-mail:libgober@math.uic.edu}

\begin{document}

\date{} 

\maketitle

\begin{abstract}
{ This paper surveys the authors' recent work on the 
two-variable elliptic genus of singular varieties. The last section
calculates a generating function for the elliptic genera of symmetric
products. This generalizes the classical results of Macdonald 
and Zagier.
}
\end{abstract}

\section{Introduction} Elliptic genera appeared in the mid-1980's 
in several diverse problems both in topology, e.g., circle actions on 
manifolds, construction of generalized
cohomology theories, genera satisfying multiplicative properties, 
and in physics, as part of the study of Dirac-like operators on loop
spaces (cf.\cite{Landweber.Stong}).
Elliptic genera are certain modular functions
attached to manifolds which interpolate many known
genera of manifolds  e.g., Todd, {\it L} and ${\it \hat A}$-genera.
Following a suggestion of E.Witten (cf. \cite{witten93}), a 
two-variable  elliptic genus
was formulated as an invariant of superconformal field theory,
and was systematically studied as a tool 
for comparison of $N=2$ minimal models and Landau-Ginzburg models 
in the work of T.Kawai,Y.Yamada and S-K. Yang (cf. \cite{KYY}).
From a mathematical point of view, the two-variable elliptic genus 
was studied in the work of Krichever, G.Hohn, B.Totaro and 
V.Gritsenko (cf. also \cite{hirz70}). 
While various generalizations were proposed (for example to complex manifolds,
cf. section \ref{start}), the two-variable elliptic
genus appears to be the most general elliptic genus in the sense that 
almost all versions of elliptic genera are its specializations.

\par The aim of these notes is to discuss generalizations
of the two-variable elliptic genus to singular varieties 
from the mathematical point of view proposed in \cite{borlibg1} and
\cite{borlibg2}, in particular without reference 
to superconformal field theories (it is curious to note,
however, the resemblance of the definition of elliptic genus in terms 
of the cohomology of the chiral de Rham complex and 
the definition of the elliptic genus of SCFT). First, we shall 
discuss the definition in terms of the cohomology of the chiral 
de Rham complex. Such a cohomology can be defined for  
hypersurfaces in Fano toric varieties in terms of the combinatorics
of the toric variety, which allows one to define the elliptic genus in 
this case. Secondly, we shall discuss the definition of elliptic genus 
of singular algebraic varieties in terms of their
resolutions and for singular spaces which are 
orbifolds $X/G$ in terms of the action of 
$G$ on $X$. These definitions can be used to 
give mathematical proofs for results which were
previously obtained from the point of view of string theory, 
notably the Dijkgraaf-Moore-Verlinde-Verlinde formula 
for the generating function of the orbifold elliptic genera
of symmetric groups acting on products of a fixed manifold $X$
(cf. section 4).
We shall finish with a derivation of generating functions for 
elliptic genera of symmetric products  
and containing as special cases old calculations of generating 
functions for Euler characteristics (I.Macdonald) and signatures
(D.Zagier).
\par This subject is extremely vast and no claim to completeness
is made. An excellent book by
F.Hirzebruch, T.Berger and R.Jung (\cite{modforms}) is 
particularly recommended for everybody interested in this 
subject.

\section{Elliptic genera of manifolds.}
\label{start}
\par Let $\Omega^{SO}_*$ (resp. $\Omega^U_*$) be the cobordism ring 
of oriented (resp. almost complex) manifolds. Recall that  
cobordism ring is defined as the quotient of the free abelian group generated by 
manifolds ($C^{\infty}$, almost complex, Spin, etc.) by the subgroup 
generated by manifolds which are 
boundaries (of manifolds with the same structure); the product is 
given by the product of manifolds. 
An $R$-valued genus is a ring homomorphism 
$E: \Omega^{SO}_* \otimes {\bf Q} \rightarrow R$.
Similarly, a complex genus is a ring homomorphism
$E: \Omega^U_* \otimes {\bf Q} \rightarrow R$.
The class of an almost complex manifold in $\Omega^{U}_* \otimes {\bf Q}$ 
is completely specified by Chern numbers (cf. \cite{Hirz1}), i.e. products of 
Chern classes evaluated on the fundamental class of the manifold.
In particular, for complex cobordism a genus can be written as 
$E(M)=\int_M {\cal E}_{{\rm dim} M}(c_1,...c_k,..)$  for some polynomial 
${\cal E}_{{\rm dim} M}$ having coefficients in the ring $R$. 
Similarly, in the oriented case, a class of $\Omega^{SO}_* \otimes {\bf Q}$ 
is determined by Pontryagin numbers and the genus is the integral 
of a polynomial in the Pontryagin classes.
\par The collection of polynomials ${\cal E}_i$  
can be specified by a characteristic series: 
$Q(x)=1+b_1x+b_2x^2+...$
($ b_i \in R$)  such that for the factorized total Chern class 
$c(T_M)=1+c_1(M)+...+c_{{\rm dim} M}(M)=(1+x_1) \cdot \cdot \cdot (1+x_r)$ 
one has ${\cal E}(c_1,...)=\prod Q(x_i)$ (cf. \cite{Hirz1}). 
This condition determines the polynomials ${\cal E}_i$ from $Q(x)$ completely. 
For example (cf. \cite{Hirz1}), the holomorphic 
Euler characteristic of a trivial bundle on a complex manifold 
extends to the complex genus and  equals the Todd genus, 
with the corresponding characteristic series being ${x \over {1-e^{-x}}}$
(Hirzebruch's Riemann-Roch theorem). 
The corresponding polynomials in Chern classes are
${c_1 \over 2}, {{c_1^2+c_2} \over 12}, {c_1c_2 \over 24}$, etc. 
In the case of oriented manifolds,  the same methods work after
replacing Chern classes by Pontryagin classes. The integer-valued 
genera which attracted the most attention, besides the Todd genus, are 
the $\it L$-genus (corresponding to the series $x \over {{\rm tanh} (x)}$; $L$-genus 
is equal to the signature of the intersection form on the middle dimensional 
cohomology cf. \cite{Hirz1}) and the ${\it\hat A}$-genus (corresponding to 
the series ${ {x /2} \over {{\rm sinh} {(x /2) }}}$ and 
equal to the index of the Dirac operator cf. \cite{AtiyahSinger}).
\par In the simplest version of the elliptic genus, the ring $R$ is the  
ring of modular forms $Mod^*(\Gamma)$ for a certain subgroup 
$\Gamma$ of 
$ SL_2({\bf Z})= \pmatrix{a \ b \cr c \ d\cr}, a,b,c,d \in {\bf Z}$,
i.e. the functions on the upper half-plane 
satisfying $\phi (\gamma \cdot \tau)=(c\tau+d)^k \phi(\tau), 
\gamma \in \Gamma$; $k$ is an integer called the {\it weight} of $\phi$
and which provides the grading of the ring of modular forms;
such functions often are written in terms of the variable 
$q=e^{2 \pi i \tau}$.  
\par Landweber-Stong (cf \cite{Landweber.Stong}) and S.Ochanine
(\cite{ochanine1}), while studying the circle actions on manifolds
and the ideals in the cobordism ring generated by the 
projectivizations of vector bundles,
considered the genus
$\Omega^* \rightarrow Mod^*(\Gamma_0(2)) \subset {\bf Q}[[q]] $,  
where $ \Gamma_0(2)= \pmatrix{a \ b \cr c \ d\cr}
\in SL(2,{\bf Z}) \vert c \ {\rm even} ) $.
Its characteristic series is given by 
\begin{equation} Q_{LSO}(x)={{x/2}\over {\sinh( x/2)}} \prod_{n=1}^{\infty}
  [{{(1-q^n)^2} \over {(1-q^n e^x)(1-q^n e^{-x})}}]^{(-1)^n}.
\end{equation}
E.Witten (\cite{witten88}) 
proposed the following expression for this genus:
$$\hat A(X) ch \{ {{R(T_X)} \over {R(1)^{{\rm dim} X}}} \}[X] $$
where 
\begin{equation} R(T_X)=\otimes_{l>0, l \equiv 0 (2)} S_{q^l}(T_X) 
   \otimes_{l>0, l \equiv 1 (2)} \Lambda_{q^l}(T_X)
\label{qbundle}
\end{equation}
and the cohomology class $ch(E)=\sum e^{x_i}$ for a bundle $E$ for which 
$c(E)=\prod(1+x_i)$ is the Chern character of $E$.
In the same paper he gave an interpretation of the elliptic genus 
as the index of a Dirac-like (or a signature-like) 
operator on the loop space ${\cal L}M$.
\par Elliptic genera of complex manifolds were defined
by F.Hirzebruch (\cite{Hirz2}) and E.Witten(\cite{witten88}).
Such an elliptic genus takes values in the ring of 
modular forms for the group
\begin{equation} \Gamma_1(N)=\{ \pmatrix{a \ b \cr c \ d\cr}
\in SL(2,{\bf Z})
\vert c \equiv 0 (N), a \equiv d \equiv 1 (N) \} 
\end{equation} 
provided the first Chern class of the manifold satisfies 
$c_1 \equiv 0 (N)$. 
\par The characteristic series depends on a choice of a point
of order $N$ on an elliptic curve with periods $2 \pi i (1,\tau)$,
say $\alpha=2 \pi i ({k \over N} \tau+ {l \over N}) \ne 0$,
 and is given in terms of 
\begin{equation} 
\Phi(x, \tau)=(1-e^{-x}) \prod_{n=1}^{\infty}
  {{(1-q^n e^x)(1-q^ne^{-x})} \over {(1-q^n)^2}}.
\end{equation}
It is equal to:
\begin{equation} 
Q_{HW}(x,\tau)=x e^{-{k \over N}x} {{\Phi(x-\alpha)} \over {\Phi(x)
     \Phi(-\alpha)}}.
\label{HW}
\end{equation}

\par I.Krichever (\cite{Krichever}) 
considered the complex genus with characteristic series
\begin{equation}
Q_K(x,z,\omega_1,\omega_2,\kappa)=
x e^{-\kappa x}
{{\sigma_{\omega_1,\omega_2}(x-z)} \over 
{\sigma_{\omega_1,\omega_2}(x)
    \sigma_{\omega_1,\omega_2}(-z)}} e^{\zeta_{\omega_1,\omega_2}(z)x}
\label{K}
\end{equation}
where $z, \kappa \in {\bf C}^*$, 
$\sigma_{\omega_1,\omega_2}(z)$ and $\zeta_{\omega_1,\omega_2}(z)$ 
are Weierstrass 
$\zeta$ ($\zeta^{\prime}=-\wp$) and $\sigma$-functions 
($\zeta={\sigma^{\prime} \over \sigma}$) 
corresponding to the same lattice 
in $\bf C$. It was further studied by G.H\"ohn (cf. \cite{hohn.thesis})
and B.Totaro (cf.\cite{Totaro}). 
In this paper B.Totaro gives an important characterization of 
the genus introduced by Krichever as the  universal genus of 
$\Omega_{SU}^*$ invariant under classical flops.
\par Note that the series $Q_K$ specializes to $Q_{HW}$ for 
 $z=\alpha$ and 
$\kappa=-{{2k} \over N} \zeta({\pi i \tau})-{{2l} \over N} 
\zeta (\pi  i)+\zeta (z)$.  In addition, the Hirzebruch-Witten genus for 
$N=2$ can be expressed in terms of Pontrjagin classes, so that 
it is  an invariant of $SO$-cobordism which  up to a factor coincides 
with the genus of Ochanine, Landweber and Stong.
\par One should mention that much of the interest in elliptic genera
first come from a conjecture by Witten later proven by
Bott and Taubes (cf. \cite{Taubes}), Hirzebruch (cf. \cite{Hirz2}), 
Krichever (cf. \cite{Krichever}) and Liu (cf. \cite{liuridg1},
\cite{liuridg}) concerning the  
{\em rigidity} property
which claims the following. Suppose a compact group $G$ acts on $M$ 
and a bundle $V$ so that an operator $P$ acting on $V$ commutes with the 
the action of $G$. Let us consider the character $L_{M,V,P}(g)=
{\rm Tr}_g {\rm Ker} P-{\rm Tr}_g {\rm Im} P$. The operator is rigid if this character is 
independent of $g$. The above mentioned results
(generalizing \cite{AtiyahHirz})  state that the 
bundles which are the coefficients of the $q$-expansion of 
(\ref{qbundle}) support operators which are rigid. 
This is the case for other genera, including (\ref{HW}) and (\ref{K}).
Another important issue in which the elliptic genus was essential is 
known under the title {\em anomaly cancellation},
which yields a series of nontrivial identities and 
congruences among various classical (i.e. $L, \hat A$ etc.)
genera (cf. \cite{liumod} and survey  \cite{Liu}). 
\par In the physics literature a two-variable elliptic genus was associated
with an $N=(2,2)$ superconformal field theory (cf.
Eguchi-Ooguri-Taormina-Yang \cite{EOTY}, E.Witten \cite{witten93} and 
Kawai-Yamada-Yang cf.\cite{KYY}). 
It is given by: 
\begin{equation} {\rm Tr}_{{\cal H}}(-1)^Fy^{J_0}q^{L_0-c/24}
\bar q ^{\bar L_0-c/24}
\end{equation}
where $\cal H$ is the Hilbert space of the SCFT, 
$L_0$ (resp. $\bar L_0$) is the Virasoro generator of left
(resp. right)-movers and $J_0$ (resp. $\bar J_0$) is the $U(1)$ charge
operator of left (resp. right) movers, the 
trace is taken over  Ramond  sector and $F=F_L-F_R$
with $F_L$ (resp. $F_R$) the fermion number of left (resp. right) movers. 
In the case when the field theory comes from a smooth 
Calabi-Yau manifold $M$, one has the following 
mathematical expression for the genus 
(cf. \cite{KYY},\cite{DMVV}, \cite{borlibg1})

\begin{equation}
Ell(M)=\int_Mch({\cal E}ll_{q,y})td(M) 
\label{charclass}
\end{equation}
where 
\begin{equation}
{\cal E}ll_{q,y}=y^{-{\dim M \over 2}}
 \otimes_{n \ge 1} (\Lambda_{-yq^{n-1}}
\bar
T_M \otimes \Lambda_{-y^{-1}q^n} T_{M} \otimes S_{q^n}\bar T_M \otimes
S_{q^n} T_M) .
\label{bundle}
\end{equation}
The characteristic series for the genus (\ref{charclass}) 
can be written in terms of the theta-function as follows.
Let  
\begin{equation}
\theta(z,\tau)=q^{1 \over 8}  (2 \sin \pi z)
\prod_{l=1}^{l=\infty}(1-q^l)
 \prod_{l=1}^{l=\infty}(1-q^l e^{2 \pi i z})(1-q^l e^{-2 \pi i
z})
\end{equation}
where $q=e^{2 \pi i  \tau}$ (the Jacobi theta-function \cite{Chandra} 
or $\theta_{1,1}$, the theta-function with theta-characteristic,
 cf. \cite{Mumford}).
Then the elliptic genus (\ref{charclass}) corresponds to the 
characteristic series (with $y=e^{2 \pi i z}$):
\begin{equation} 
x \cdot  {{\theta ({{x} \over {2 \pi i}}-z,\tau)} \over
{\theta ({{x} \over {2 \pi i }}, \tau)}}
\label{charseries}
\end{equation}
(cf. \cite{KYY} and \cite{borlibg1}). 
Note that the use of theta-functions in connection with 
elliptic genera goes back to D.Zagier (cf. \cite{zagier})  
and J.L. Brylinski (\cite{bryl}). 
\par The elliptic genus $K(M,\omega_1,\omega_2,z,\kappa)$ introduced by 
I.Krichever for a Calabi-Yau manifold $M$ differs from 
the elliptic genus (\ref{charclass}) only by a factor 
which depends only on dimension (and is independent of $\kappa$ 
cf. \cite{borlibg1} Sect.2):  
\begin{equation} K(2\pi \ii z, 2 \pi \ii, 2 \pi \ii \tau,
\kappa)(X)=
Ell(z,\tau)(X) \cdot (-{{\theta^{\prime}(0,\tau)} 
\over {2\pi \ii\, \theta(z,\tau)}})^d.
\end{equation} 

\par The automorphic property of the elliptic genus is central for 
understanding this invariant. Recall that a weak Jacobi form of weight $k$ 
and index $r$ ($k \in {\bf Z}, r \in {1 \over 2} {\bf Z}$: we consider 
forms of half-integral index) is a holomorphic function on $H \times {\bf C}$ 
satisfying: 
\begin{equation} \phi({{a\tau +b} \over {c \tau +d}},{z \over {c \tau +d}})=
 (c \tau +d)^k e^{2 \pi i{{r cz^2} \over {c \tau +d}}  }\phi (\tau, z)
\end{equation}
\begin{equation} \phi(\tau, z+m\tau+n)=
(-1)^{2r(\lambda+\mu)} e^{-2 \pi i r (m^2 \tau +2 m z)} \phi (\tau, z)
\end{equation}
In addition, a weak Jacobi form must have a Fourier expansion with non-negative
powers of $q=e^{2\pi i \tau}$. This is weaker than the usual condition on
Fourier modes, which explains the name (cf. \cite{EZ}).
\par Using the expression via $\theta$-functions for the characteristic series of
the elliptic genus (\ref{charseries}), one can show that 
the elliptic genus of an (almost) complex manifold of dimension $d$ 
is a weak Jacobi form of weight $0$, 
index $d \over 2 $ (cf. \cite{borlibg1}).
A description of the space of weak Jacobi forms in \cite{EZ} yields that 
elliptic genera of Calabi-Yau manifolds span the space 
of Jacobi forms of weight 0 and index $d \over 2$ (cf. \cite{borlibg1}, theorem 2.6).
Gritsenko (\cite{Gritsenko}) has calculated the $\bf Z$-span of elliptic genera. 
\par Such calculations in particular allow one to decide to what extent 
the elliptic genus depends on the $\chi_y$-genus. Note that the elliptic genus 
is a combination of Chern numbers and there are non-trivial relations 
among Chern and Hodge numbers (e.g. 
$\sum_{p=2}^d (-1)^p {p \choose 2} \chi^p={1 \over 12}\{{1 \over2 }
d(3d-5)c_d+c_{d-1}c_1 \}[X]$ cf. \cite{LW}). 
More precisely,
\begin{theo} If the dimension of a Calabi-Yau manifold is less than 12 or is
equal to 13, then the numbers $\chi_p$ determine its elliptic genus
uniquely. In all other dimensions there exist Calabi-Yau manifolds with
the same $\{\chi_p\}$ but distinct elliptic genera. 
\label{hodge}
\end{theo}
For example, if $e(X)$ (resp. $\chi(X)$) denotes the topological 
(resp. holomorphic) Euler characteristic
then (cf. \cite{Neumann}, \cite{KYY})
the elliptic genus in the case of threefolds is
\begin{equation}
{e(X) \over 2}(y^{-{1 \over 2}}+y^{1 \over 2})\prod_{n=1}^{n=\infty}
{{(1-q^ny^2)(1-q^ny^{-2})} \over {(1-q^ny)(1-q^ny^{-1})}}
\end{equation}
and for fourfolds is
\begin{equation}
\chi(X) E_4A^2+{{e(X)} \over 144}(B^2-E_4A^2).
\end{equation}
Here, $A={{\phi_{10,1}(\tau,z)} \over {\eta^{24}(\tau)}}, 
B={{\phi_{10,1}(\tau,z)} \over {\eta^{24}(\tau)}}$, where
$\phi_{10,1}$ and $\phi_{12,1}$ are the unique cusp forms of index 1 and
weights 10 and 12, resp. (cf.\cite{EZ}), $\eta(\tau)$ is the Dedekind 
$\eta$-function and $E_4(\tau)$ is the normalized Eisenstein series of 
weight 4.
\par \noindent 
However, as follows from the above theorem,
for manifolds of high dimension
the elliptic genus contains information not available from the $\chi_y$-genus.
It is interesting, therefore, to know 
what are the values of this invariant for concrete manifolds. For example,  
the $\chi_y$ characteristic of toric varieties is well known (cf. \cite{Danilov},
\cite{Oda} or \cite{Fulton}). For elliptic genera of smooth toric varieties we
have the following:
\begin{theo} Let $\bf P$ be a smooth toric variety corresponding to a fan 
$\Sigma$ in $N \otimes {\bf R}$ for some lattice of rank $d$. Let $M$ be the 
lattice dual to $N$. For the cone  $C^*$ 
of $\Sigma$ (which is simplicial due to the smoothness of $\bf P$)
let $n_i (i=1,...,d)$ be a system of its  generators.
Then:   
\begin{equation}Ell({\bf P},y,q)=y^{-d/2}
\sum_{m\in M}
\sum_{C^*\in\Sigma}
(-1)^{{\rm codim} C^*}
\left(
\prod_{i=1,...,\dim C^*} {1\over 1-yq^{m\cdot n_i}}
\right)
G(y,q)^d
\label{toricgenus}
\end{equation}
where
$$G(y,q)=\prod_{k\geq 1}
{(1-yq^{k-1})(1-y^{-1}q^k)\over(1-q^k)^2}.$$
\end{theo}

We shall  sketch the proof, which uses the calculation of the cohomology 
via a split of the \v{C}ech complex according to characters.
\par First, let us consider the Leray spectral sequence for 
the cover of the toric variety by open sets
 ${\bf A}_C={\rm Spec} {\bf C}[C]$ defined by the cones $C^* \in \Sigma$
and apply this spectral sequence to the bundle ${\cal E}ll_{q,y}({\bf P})$ 
 (cf.(\ref{bundle})).
By abuse of language, the bundle here actually is a bigraded  bundle
whose components are the coefficients of $y^aq^b$ in 
${\cal E}ll_{q,y}({\bf P})$;
these coefficients are bundles having finite rank.
Since the cohomology of positive dimension of the bundle 
${\cal E}ll_{q,y}({\bf P})$ 
vanishes over affine sets, it yields 
$$Ell({\bf P};y,q)=y^{-d/2}\sum_{m\in M}
(\sum_{C_0^*,...,C_k^*}
(-1)^{k} \dim_m H^0({\bf A}_{C_0}\cap ... \cap {\bf A}_{C_k},{\cal E}ll_{q,y}
({\bf P})).
$$
Second, over each such open set ${\bf A}_C$ of maximal dimension, 
since  ${\bf A}_C$ is just 
an affine space, a direct calculation shows 
\begin{equation} \sum_{m\in M} 
t^m \dim_m H^0({\bf A}_C,{\cal E}ll_{q,y}({\bf P}))=
\prod_{i=1,...,d} \prod_{k\geq 1} 
{(1-t^{m_i}yq^{k-1})(1-t^{-m_i}y^{-1}q^k)
\over
(1-t^{m_i}q^{k-1})(1-t^{-m_i}q^k).
}
\label{toric1}
\end{equation} 
where $m_i$ are generators of the cone $C$ forming a basis in the
lattice of the space containing the cone.
Third, one notices that the latter can be rewritten as
\begin{equation}
\prod_{i=1,...,d} \prod_{k\geq 1}
{(1-t^{m_i}yq^{k-1})(1-t^{-m_i}y^{-1}q^k)
\over
(1-t^{m_i}q^{k-1})(1-t^{-m_i}q^k)
}
= \sum_{m\in M} t^m
\prod_{i=1,...,d} \left({1\over 1-yq^{m\cdot n_i}}\right)
G(y,q)^d
\label{toric2}
\end{equation}
where
$$G(y,q)=\prod_{k\geq 1}
{(1-yq^{k-1})(1-y^{-1}q^k)\over(1-q^k)^2}$$
and $n_i$ are generators of $C^*$. 

One checks that the combined result of (\ref{toric1}) and (\ref{toric2})
is true for cones of arbitrary (i.e. possibly nonmaximal) dimension
since ${\bf A}_{C}={\bf C}^{{\rm dim}C^*} \times ({\bf C}-0)^{d-{\rm dim}C^*}$.

\par Finally, a combinatorial argument shows that the total 
contribution of each cone in the terms of the \v Cech complex,
i.e. $\sum_{C_0 \cap ... \cap C_k=C}  (-1)^{k}$,
 is equal to $(-1)^{{\rm codim} C^*}$. This yields the 
theorem. 
\par Since compact toric varieties are never Calabi-Yau, the expression 
(\ref{toricgenus}) is not expected to have automorphic properties.
However, its specialization to one-variable genera must satisfy 
modular relations.
For example, for the Landweber-Stong-Ochanine elliptic genus
  $$\widehat{Ell}(X;q)=(-1)^{d/2} 
   Ell(X;-1,q) G(-1,q)^{-d}$$
we obtain 
\begin{theo}
If $\bf P$ is a smooth complete toric variety, then 
$$\widehat{Ell}({\bf P};q)=\sum_{m\in M} \left(
\sum_{C^*\in\Sigma}
(-1)^{{\rm codim}\  C^*}
\prod_{i=1,...,\dim \  C^*} {1\over 1+q^{m\cdot n_i}}
\right).
$$
In particular, the series in the right hand side is a modular form.
\end{theo}

It is interesting that neither the modular property nor  
the relation to previous calculations of 
elliptic genera are obvious but, rather, lead to interesting 
new identities. For example, since
$${\widehat {Ell}}({\bf P}^2)=\delta=-{1 \over 8}-3 \sum_{n \ge 1}
(\sum_{d \vert n, d \ {\rm odd}} d)q^d$$ 
we have
$$\sum_{m\geq 1,n\geq 1} 
{q^{m+n}\over (1+q^m)(1+q^n)(1+q^{m+n})}=
\sum_{r\geq 1}q^{2r}\sum_{k|r}k. $$
(cf. \cite{borlibg1} for a direct proof of this identity, rather 
than as a consequence of two different calculations of elliptic genera).
\par The next problem is how to calculate the elliptic genus of hypersurfaces 
in toric varieties. To describe this, one needs a
description of the elliptic genus via the  
chiral de Rham complex.

\section{Elliptic genera in the singular case and the chiral de Rham complex}

The two-variable elliptic genus is closely related to the chiral de 
Rham complex
constructed by Malikov, Schechtman and Vaintrob in \cite{MSV} for algebraic 
(analytic, $C^{\infty}$ etc.) manifolds. This is a sheaf 
of vector spaces which has the structure of sheaf of vertex operator
algebras. In particular, it supports the action of the Virasoro
algebra, whose role in the theory of elliptic genera was anticipated from the
very beginning  (cf. \cite{segal}; for another attempt to clarify the 
role of the Virasoro algebra cf. \cite{tomanoi}). 
\par For convenience, let us recall the definition of a vertex operator algebra
and conformal vertex operator algebra (cf. for example \cite{Kac}).
\begin{dfn}A vertex operator algebra is a vector space $V$, endowed with 
\par \noindent 1. a decomposition 
\begin{equation} V=V_0 \oplus V_1
\label{split}
\end{equation}
\par \noindent 2. a vector denoted $ \vert 0> \in V_0$
and called the vacuum vector 
\par \noindent 3. a linear map $V \rightarrow End(V)[z,z^{-1}]$ called the 
states to 
fields correspondence; the image of $a \in V$ is denoted 
$Y(a,z)=\sum_{n \in {\bf Z}} a_{(n)}z^{-n-1}, a_{(n)} \in End(V)$. One
requires that for fixed $a$ and $b$ there holds $a_{(n)}b=0$ for $n \gg 0$.
\par \noindent 4. a linear map $T:V \rightarrow V$ called the 
infinitesimal translation 
operator. 
\par This data are required to satisfy the following axioms:
\par \noindent a)Translation covariance: $\{T,Y(a,z)\}_{-}=\partial Y(a,z)$.
\par \noindent b)Vacuum: $ \vert 0>$ satisfies: $Y(  \vert 0 >,z)=
{\rm Id}_V$,
$Y(a,z)  \vert 0 > \vert_{z=0} =a$, $T \vert 0>=0$

\par \noindent c)Locality: $(z-w)^NY(a,z)Y(b,z)=(-1)^{p(a)p(b)}(z-w)^N
Y(b,z)Y(a,z)$ for $N \gg 0$
\end{dfn}
\begin{dfn}
\par A conformal vertex algebra is a pair $(V,L)$, where $V$ is a  
vertex algebra and $L$ is a field that corresponds to an even element
with the following properties:  
\par \noindent 1. The components of $L(z)=\sum_nL_nz^{-n-2}$
satisfy the Virasoro commutation relations:
$$[L_n,L_m]=(n-m)L_{n+m}+{{n^3-n} \over 12} \cdot c \cdot \delta^n_{-m}$$ 
\par \noindent 2. $L_{-1}=T$ is an infinitesimal translation operator.
\par \noindent 3. $L_0$ is diagonalizable. 
\end{dfn}

In \cite{MSV} the authors prove the following: 
\begin{theo} Let $X$ be a nonsingular compact complex manifold. 
There exists a sheaf $\Omega^{ch}_X$ of vector spaces on $X$ with the 
properties: 
\par a) For each Zariski open set $U$, $\Gamma(U,\Omega^{ch}_X)$ 
has a structure of conformal vertex algebra, with the restriction maps
being morphisms of vertex algebras. 
\par b) $\Omega^{ch}_X$ has two gradings with degrees called
 fermionic charge and conformal weight.
\par c) $\Omega^{ch}_X$ has de Rham differential $d_{DR}^{ch}$ of 
(fermionic) degree 1,  $(d_{DR}^{ch})^2=0$.

\par d) The usual de Rham complex $\Omega_X^{\bullet}$ is isomorphic 
to the conformal weight zero component of   $\Omega_{DR}^{ch}$.

\par e)The complex $(\Omega_X^{ch}, d_{DR}^{ch})$ is 
quasi-isomorphic to $(\Omega_X^{\bullet},d_{DR})$.

\par f) Each component of fixed conformal weight 
has a canonical filtration with $gr_F$ isomorphic to 
the tensor product of exterior powers of the tangent and cotangent bundles,
so that the corresponding generating function is 
 $$\otimes_{n \ge 1} (\Lambda_{yq^{n-1}}
\bar T_X \otimes \Lambda_{y^{-1}q^n} T_X \otimes S_{q^n}\bar T_X
\otimes S_{q^n}T_X)$$
\end{theo}
Recall that the supertrace of an operator $S$ acting on a space (\ref{split}) 
is ${\rm tr} S \vert_{V_0}-{\rm tr} S \vert_{V_1}$. 
By the Riemann-Roch theorem, the integral in (\ref{charclass}) is just 
$\sum_i (-1)^i {\rm dim} H^i({\cal E}ll_{q,y}(M))$. 
If one considers the bigraded 
sheaf with components being the coefficients of (\ref{charclass}), the parity 
given by the parity of the exponent of $y$  and endowed 
with the operators $A$ and $B$ acting on the coefficient of of $y^aq^b$ as 
multiplication by $a$ and $b$ respectively, then     
we see that elliptic genus can be written as 
$y^{{-{\rm dim} M} \over  2}
{\rm Supertrace}_{H^*({\cal E}ll_{q,y}(M))}y^Aq^B$. 
Since the Euler characteristics of a filtered sheaf and its associated graded 
sheaf are the same, this suggests the following: 
\begin{dfn} Let $X$ be a variety for which one can define 
a chiral de Rham complex $\Omega^{ch}_X={\cal MSV}(X)$ with properties a)-f) 
as above.
The elliptic genus of $X$ is then defined as
 $$y^{-{{{\rm dim} X} \over 2}}{\rm SuperTrace}_{H^*({\cal MSV}(X))}y^{J[0]}q^{L[0]}.$$
\end{dfn}

The usefulness of this definition stems from the following:
the first-named author did construct such a complex ${\cal MSV}(X)$ in the 
case when $X$ is a hypersurface in a toric varieties
with Gorenstein singularities (cf. \cite{Bvertex}) or for 
toric varieties themselves. In \cite{Bvertex}, a purely combinatorial
construction of the cohomology of ${\cal MSV}(X)$ is given in these 
cases. It contains a description of the latter as 
the BRST cohomology of Fock spaces with an explicit description of those
in terms of combinatorics. This yields the following explicit 
formulas for elliptic genera.
\begin{theo}
{\rm
Let $X$ be a generic hypersurface in a Gorenstein toric Fano variety
corresponding to a reflexive polytope $\Delta$ in a lattice $M_1$,
${\rm rk}M_1= d+1$.
Let $M=M_1 \oplus {\bf Z}$,  $N_1$ and 
$N$ be the lattices dual to $M_1$ and $M$ respectively and $\Delta^*$
be the polytope dual to $\Delta$. Denote 
the elements $(0,1) \in M, (0,1) \in N$ as $\rm deg$ and $\rm deg^*$ 
respectively. Let $K$ (resp. $K^*$) be the cone in $M$ (resp. $N$)
over $(\Delta,1)$ (resp. $(\Delta^*,1)$) with the vertex at 
$(0,0)_M$ (resp. $(0,0)_N$). 
Then
$$
Ell(X,y,q)=
y^{-{d \over 2} } \sum_{m\in M} 
\left(
\sum_{n\in K^*}
y^{n\cdot {\deg} -m\cdot {\deg}^*}
q^{m\cdot n+m\cdot{\deg}^*}
G(y,q)^{d+2}
\right)
$$
where 
$$G(y,q)=\prod_{k\geq 1}
{(1-yq^{k-1})(1-y^{-1}q^k)\over(1-q^k)^2}.$$}
\label{ellgenhyp}
\end{theo}

\noindent On the other hand, in the toric case one obtains: 

\begin{theo} {\rm For a toric Gorenstein variety ${\bf P}$ 
$$
Ell({\bf P},y,q)=y^{-d/2}
\sum_{m\in M}
\sum_{C^*\in\Sigma}
(-1)^{{\rm codim} C^*} 
(\sum_{n\in C^*}q^{m\cdot n}y^{\deg\cdot n})
G(y,q)^d.
$$}
\end{theo}

\noindent The Gorenstein property is needed since ${\bf P}$ has 
Gorenstein singularities if and only if 
the function of $n$ given by $n \cdot {\rm deg}$ takes integer values. Inspection of 
the formulas in these theorems yields the following: 

\begin{coro} {\rm If $X$ admits a crepant toric desingularization 
 $\hat X$, then
 $$Ell(X,y,q)=Ell(\hat X,y,q).$$}
\end{coro}

Similarly to the non singular case we have:
\begin{theo}{\rm 
 The elliptic genus of a generic Calabi-Yau hypersurface
in a any toric Gorenstein Fano variety is a weak Jacobi form of 
weight 0 and index $d \over 2 $.}
\end{theo}

\noindent The proof uses an extension of the elliptic genus to a 
three-variable function and expression of the latter via theta 
functions, which reduces to the Bott formula in the smooth case.
(cf. lemma 5.3 in \cite{borlibg1})

\par \noindent An explicit description of the Fock spaces for which 
BRST cohomology yields the cohomology of the chiral 
de Rham complex ${\cal MSV}(X)$ in the case of theorem \ref{ellgenhyp}
and use of the Jacobi property of their elliptic genus 
provides the following relation: 

\begin{theo}Let {\rm $X,X^*$ be Calabi-Yau hypersurfaces in toric Gorenstein
Fano varieties corresponding to dual reflexive polytopes
$\Delta$ and $\Delta^*$. Then:
$$Ell(X;y,q)=(-1)^d Ell(X^*;y,q).$$}
\label{mirror}
\end{theo}

\par \noindent Such a result certainly is expected from physics 
considerations and assuming that Calabi-Yau hypersurfaces 
corresponding to dual polytopes form a mirror pair in the
strong sense of correspondence between CFT's.
Also, one can check it in small dimensions when the elliptic 
genus is a combination of Hodge numbers (cf. \ref{hodge}
and \cite{Neumann} for explicit formulas). But in higher
dimensions the relation in theorem \ref{mirror} can be viewed
as a test for deciding if two Calabi-Yau manifolds
form a mirror pair.

\section{Elliptic genus of singular varieties via resolution 
of singularities and orbifold elliptic genera.}

The definition of elliptic genus for special singular varieties
in the last section suggests the following problem: 
find an expression for the elliptic
genus of singular varieties  
in terms of a resolution and define the elliptic genus for 
varieties more general than hypersurfaces in singular toric  spaces.
These problems were  addressed in \cite{borlibg2}, where the
following approach was proposed.

\begin{dfn} {\rm Let $Z$ be a complex space with $\bf Q$-Gorenstein 
singularities and let $Y \rightarrow Z$ be a resolution of 
singularities. Let $\alpha_k \in {\bf Q}$ be the discrepancies, i.e.
rational numbers defined from the relation:
$K_{Y}=\pi^*K_{Z}+\sum \alpha_kE_k$. Then 
$$
{Ell}_{sing}(Z;z,\tau):=
\int_Y
\Bigl(\prod_l
\frac{(\frac{y_l}{2\pi\ii})\theta(\frac{y_l}{2\pi\ii}-z,\tau)\theta'(0,\tau)}
{\theta(-z,\tau)\theta(\frac{y_l}{2\pi\ii},\tau)}
\Bigr)\times
\Bigl(\prod_k
\frac{\theta(\frac{e_k}{2\pi\ii}-(\alpha_k+1)z,\tau)\theta(-z,\tau)}
{\theta(\frac{e_k}{2\pi\ii}-z,\tau)\theta(-(\alpha_k+1)z,\tau)}
\Bigr)
$$ }
\label{sing}
\end{dfn}

\noindent (This definition can be generalized
to define the elliptic genus of log-terminal pairs; cf. \cite{borlibg2}
for details).

\par It turns out that $Ell_{sing}(Z;z,\tau)$ is independent of $Y$
and hence defines an invariant of $Z$. Several results make 
this invariant interesting. 
\bigskip 
\par \noindent 1. It does specialize to the 
normalized version of the elliptic genus discussed earlier in the case 
when $Z$ is non-singular, i.e. 
\begin{equation}Ell_{sing}(Z,z,\tau)=
Ell(Z,z,\tau) (-{{\theta^{\prime}(0,\tau)} 
\over {2\pi \ii\, \theta(z,\tau)}})^d
\label{normalizing}
\end{equation}  
\bigskip \par  \noindent 2. If $Z$ admits a crepant resolution,
 i.e. such that all discrepancies
are zero, the singular elliptic genus coincides with the 
elliptic genus of a crepant resolution (up to the same factor as in 
(\ref{normalizing})).    

\bigskip \par \noindent 3. In the case when $Z$ is a Calabi-Yau,
the singular elliptic genus has the transformation properties of a 
Jacobi form. 
\bigskip 
\par \noindent 4. For Calabi-Yau hypersurfaces in Fano Gorenstein toric 
varieties, the elliptic genus in \ref{sing} coincides with the elliptic 
genus considered in the last section 
(again up to the factor in (\ref{normalizing})).
\bigskip 
\par \noindent  5. If $q \rightarrow 0$ then the singular elliptic genus 
specializes (up to a factor) to the $\chi_y$ genus that is a 
specialization of the
$E$-function  studied by Batyrev (\cite{Batyrev}).

Finally, in many situations $Ell_{sing}$ is related to the
elliptic genus of orbifolds, also introduced in \cite{borlibg2}.
Let $X$ be a complex manifold on which a finite group $G$
is acting via holomorphic transformations. 
Let $X^h$ will be the fixed point set of $h \in G$
and  $X^{g,h}=X^g \cap X^h (g,h \in G)$. Let 
\begin{equation}
TX \vert_{X^h}= \oplus_{\lambda(h) \in \QQ \cap [0,1)} V_{\lambda}.
\label{directsum}
\end{equation}
where the bundle $V_{\lambda}$ on $X^h$ is determined by the 
requirement that $h$ acts on $V_{\lambda}$ via multiplication by 
$e^{2 \pi i \lambda(h)}$. For a connected component of $X^h$
(which by abuse of notation we also will denote $X^h$), the fermionic 
shift is defined as $F(h,X^h \subseteq X)=\sum_{\lambda} \lambda(h)$
(cf. \cite{zaslow}, \cite{Batyrev.Dais}). Let us consider the 
bundle:
$$
V_{h,X^h\subseteq X}:= \otimes_{k\geq 1} 
\Bigr[
(\Lambda^\bullet V_0^*yq^{k-1})\otimes
(\Lambda^\bullet V_0  y^{-1}q^{k})\otimes
(Sym^\bullet V_0^*q^{k})\otimes
(Sym^\bullet V_0 q^{k})\otimes
$$
\begin{equation}\otimes
\bigl[
\otimes_{\lambda\neq 0} 
(\Lambda^\bullet V_\lambda^*yq^{k-1+\lambda(h)})\otimes
(\Lambda^\bullet V_\lambda  y^{-1}q^{k-\lambda(h)})\otimes
(Sym^\bullet V_\lambda^*q^{k-1+\lambda(h)})\otimes
(Sym^\bullet V_\lambda q^{k-\lambda(h)})
\bigr]
\Bigl]
\end{equation}
\begin{dfn} {\rm The orbifold elliptic genus of a $G$-manifold $X$
is the function on $H \times {\bf C}$ given by:
$$Ell_{orb} (X,G;y,q) := 
y^{-\dim / X /2} \sum_{\{h\},X^h}   y^{F(h,X^h\subseteq X)} \frac 1 {|C(h)|}
\sum_{g\in C(h)} L(g, V_{h,X^h\subseteq X})
$$
where the summation in the first sum is over all conjugacy classes in 
$G$ and connected components $X^h$ of an element $h \in \{h\}$,
$C(h)$ is the centralizer of $h \in G$ and 
$L(g,V_{h,X^h\subseteq X})=
\sum_i (-1)^i{\rm tr}(g, H^i(V_{h,X^h\subseteq X}))$ is the holomorphic
Lefschetz number.}
\end{dfn}

Using the  holomorphic Lefschetz formula 
(\cite{AtiyahSinger}) one can rewrite
this definition as follows.

\begin{theo} {\rm Let $TX \vert_{X^{g,f}}= \oplus W_{\lambda}$ and let  
$x_{\lambda}$  be the collection of Chern roots of $W_{\lambda}$. 
Let $$\Phi(g,h,\lambda,z,\tau, x)=
{{\theta({x \over 2 \pi i }+\lambda(g)-\tau \lambda(h)-z)}
 \over 
{\theta({x\over 2 \pi i }+\lambda(g)-\tau \lambda(h))}}
e^{2\pi i z \lambda(h)}.$$
Then:$$E_{orb}(X,G,z,\tau)={1 \over {\vert G \vert}} 
\sum_{gh=hg} \prod_{\lambda(g)=\lambda(h)=0}
x_{\lambda}\prod_{\lambda} \Phi(g,h,\lambda,z,\tau,x_\lambda)[X^{g,h}]. $$
}
\label{commute}
\end{theo}

An orbifold elliptic genus so defined specializes for $q=0,y=-1$ to 
the orbifold Euler characteristic: 
$ e_{orb}(X,G)={1 \over {\vert G \vert}}\sum_{fg=gf} e(X^{f,g})$ 
(cf. \cite{hirzhofer} and \cite{atiyahsegal} where such 
an orbifold Euler characteristic is interpreted as the Euler 
characteristic of equivariant $K$-theory: ${\rm rk} K^0_G(X)-
{\rm rk}K^1_G(X)$).
Such an orbifold elliptic genus also can be specialized to the 
orbifold $E$-function studied by Batyrev-Dais (\cite{Batyrev.Dais}).
Moreover, one can show that $Ell_{orb}(X,G)$ is an invariant of cobordisms
of $G$-actions.
\par One of the consequences of \ref{commute} is the Jacobi property
of $Ell_{orb}(X,G)$ in the case when $X$ is Calabi-Yau and 
the action of $G$  preserves a holomorphic volume form 
(for more general actions, $Ell_{orb}$ still has the Jacobi 
property but only for a subgroup of the Jacobi group described 
in terms of the order of the image of $G$ in${\rm Aut} H^0(X,{\Omega^d(X)})$).
\par In the case when $X \rightarrow X/G$ does not have ramification 
we have the following conjecture. 
\begin{conj}\label{mainconj}
{\rm
Let $X$ be a complex manifold equipped with an effective action 
of a finite group $G$. Then
$$
Ell_{orb}(X,G;y,q)=
\left(
\frac 
{2\pi\ii\theta(-z,\tau)}
{\theta'(0,\tau)}
\right)^{\dim X}
\widehat{Ell}(X/G,;y,q)
$$
}
\end{conj}
(for a more general statement, which allows ramification,
cf. \cite{borlibg2}). This conjecture is proven in \cite{borlibg2}
in the case when $X$ is a smooth toric variety and $G$ is 
a subgroup of the big torus and also for arbitrary $X$ in the case 
when $G={\bf Z}/2{\bf Z}$
(using the description of generators of the cobordisms of 
${\bf Z}/2{\bf Z}$-actions given in \cite{kosn}). 
Assuming this conjecture in the case when  
$X/G$ admits a crepant resolution $\widetilde {X/G}$, the orbifold ellptic
genus is just the elliptic genus of such a resolution. 
So it is natural to think about $Ell_{orb}(X,G)$ as 
a substitute for the elliptic genus of crepant resolution in the 
cases when it does not exist.  
\par The most interesting property of $E_{orb}(X,G)$ is that 
it yields the remarkable formula due to R. Dijkgraaf, G. Moore, E. Verlinde 
and H. Verlinde (cf. \cite{DMVV}, also cf. \cite{dijkgraaf})
obtained as part of the identification of the elliptic 
genus of the supersymmetric sigma model of the $N$-symmetric product
of a manifold $X$ and
the partition function of a second quantized string theory on 
$X \times S^1$. Namely, in \cite{borlibg2} a mathematical 
proof is given for the following.

 \begin{theo}
Let $X$ be a smooth variety $X$ with elliptic genus
$\sum_{m,l}c(m,l)y^lq^m$. Then 
$$
\sum_{n\geq 0}p^n Ell_{orb}(X^n,\Sigma_n;y,q)=
\prod_{i=1}^\infty \frac 1
{(1-p^iy^lq^m)^{c(mi,l)}}.
$$
\label{DMVVtheorem}
\end{theo}

Note that since the elliptic genus can be specialized to $\chi_y$-genus
and the Hilbert schemes for surfaces give a crepant resolution of 
the symmetric product, the results of \cite{gottell} and 
\cite{gottsche} can be viewed as special cases of this theorem
(cf. also \cite{zhou2}).

\section{Generating functions for elliptic genera of symmetric products.}
Another interesting question is about a generating function 
similar to \ref{DMVVtheorem} 
but constructed
for ordinary elliptic genus of the quotient which we define as
\begin{equation}
{1 \over {\vert G \vert }} \sum_{g} L(g,{\cal E}ll_{q,y}(X))
\end{equation}
where ${\cal E}ll_{q,y}(X)$ is the bundle (\ref{bundle}).
We remark that this represents a ``naive'' version of an elliptic
genus of the quotient, and is {\it different} from the orbifold genus 
considered in the last section. In particular, one cannot expect it to 
satisfy the formula of \cite{DMVV}. 
On the other hand, such an elliptic genus of the quotient specializes to 
the $\chi_y$-genus of the quotient (cf. \cite{zagier}) 
and in particular determines the Euler characteristic and
the signature of the quotient.  
Generating functions for these classical invariants 
of symmetric products of manifolds 
were obtained earlier:   
for the Euler characteristic (cf. (\ref{macdnfla}) and \cite{macdonald})
and for the signature (cf. \cite{zagier}, \cite{zhou1},\cite{zhou2},
and \ref{zagierfla}).
The analog of  \ref{DMVVtheorem}
is the following: 
\begin{theo}
Let $Ell(X)=\Sigma c(m,l)q^my^l$. Then 
$$\sum_n
 Ell(X^n/\Sigma_n)t^n=
\prod_{m,l}{1 \over {(1-tq^my^l)^{c(m,l)}}}.
$$
\label{generating}
\end{theo} The proof is based on the following expression of holomorphic
Lefschetz numbers of (\ref{bundle}) via theta functions.

\begin{lem}
$$L(g,y^{-d/2} \Lambda_{-yq^{k-1}}T^* \otimes
 \Lambda_{-y^{-1}q^k}T \otimes S_{q^k} (T^*) \otimes S_{q^k}(T))$$
$$= {\prod_{i,r,s} {{y_i \theta({{y_i} \over {2 \pi i }}-z,\tau) 
\theta ({{{x_{r,s}+\theta_r}} \over
{2 \pi i }}-z,\tau)}} \over 
{\prod_{r,s,i} \theta ({{y_i} \over {2 \pi i }},\tau) 
\theta ({{x_{r,s}+\theta_r} 
\over {2 \pi i}},\tau)}}$$

\end{lem}

\par \noindent {\bf Proof.}  We shall use the Atiyah-Singer holomorphic 
Lefschetz formula:
 $$L(g,V)=
{{[ch \ V \vert_{X^g}](g) td(T_{X^g}) } \over {ch \ \lambda_{-1} 
 (N^g)^*(g)}}[X^g]$$
If $N^g=\oplus N^g(\theta_r)$ has Chern roots $x_{r,s}$ then 
 $ch \ \lambda_{-1}((N^g)^*)(g)=\prod_{r,s} (1-e^{-x_{r,s}-\theta_r})$ 
Let $y_i$ be the Chern roots of $T_{X^g}$.
Then we have: 
$${{ch {\cal E}ll_{q,y}(X) \vert_{X^g} td(X^g)} \over {ch \lambda_{-1}((N^g)^*)(g)}}
=$$ $$y^{-d/2}{{\prod_{i,r,s} y_i(1-yq^{k-1}e^{-y_i})
(1-yq^{k-1}e^{-x_{r,s}- \theta_r})
(1-y^{-1}q^{k}e^{y_i})(1-y^{-1}q^ke^{x_{r,s}+\theta_r})} \over 
{\prod_{i,r,s} (1-q^ke^{-y_i})(1-q^ke^{-x_{r,s}-\theta_r})(1-q^ke^{y_i})
(1-q^ke^{x_{r,s}+\theta_{r,s}})(1-e^{-y_i}) 
\prod_{r,s} (1-e^{-x_{r,s}-\theta_r})} }$$ 
The latter can be written as
$$y^{-d/2}{{\prod_{i,r,s} y_i(1-yq^{k}e^{-y_i})
(1-yq^{k}e^{-x_{r,s}- \theta_r})
(1-y^{-1}q^{k}e^{y_i})(1-y^{-1}q^ke^{x_{r,s}+\theta_r})(1-ye^{-y_i})
 (1-ye^{-x_{r,s}-\theta_r} )} \over 
{\prod_{i,r,s} (1-q^ke^{-y_i})(1-q^ke^{-x_{r,s}-\theta_r})(1-q^ke^{y_i})
(1-q^ke^{x_{r,s}+\theta_{r,s}})(1-e^{-y_i}) 
\prod_{r,s} (1-e^{-x_{r,s}-\theta_r})}} $$
Since $\sin\pi (a-z)=e^{\pi i (a-z)} {(1-e^{-2 \pi i (a-z)}) \over {2
i}}=y^{-{1 \over 2}} e^{\pi i a}(1-ye^{-2 \pi i a})({1 \over {2i}})$ 
this can be written as:
$$\prod_{i,r,s} {{2 \sin \pi ({{y_i} \over {2 \pi i}}-z)(1-e^{2 \pi i z}q^ke^{-y_i})
(1-e^{-2 \pi i z} q^k e^{y_i})2 \sin\pi ({{{x_{r,s}+\theta_r}} \over {2 \pi i}}-z) 
(1-e^{2 \pi i z}q^k
e^{-x_{r,s}+\theta_r+2 \pi i z}}) \over 
{2 \sin \pi {y_i} (1-q^ke^{y_i})(1-q^ke^{-y_i})2 \sin \pi (x_{r,s}+\theta_r) 
(1-q^ke^{x_{r,s}+\theta_r})(1-q^ke^{-x_{r,s}-\theta_{r,s}})}}
 $$
$$ {{(1-e^{2 \pi i z}q^k e^{-x_{r,s}+\theta_r})} \over }=
 {{{\prod_{i,r,s} {{y_i \theta({{y_i} \over {2 \pi i }}-z,\tau) 
\theta ({{{x_{r,s}+\theta_r}} \over
{2 \pi i }}-z,\tau)}}} \over { 
{\prod_{r,s,i} \theta ({{y_i} \over {2 \pi i }},\tau) 
\theta ({{x_{r,s}+\theta_r} 
\over {2 \pi i}},\tau)}}}}.$$

\hfill{$\Box$}

We also shall use the following two identities:

$$\prod_{k=0}^{k=r-1} \sin \pi (x+{k \over r})= 
{1 \over {2^{r-1}}}\sin \pi r x $$
and $$\prod_{k=0}^{k={r-1}} (1-q^le^{2 \pi i z +2 \pi i {k \over r}})=
(1-q^{rl}e^{2 \pi i z r})$$ 
which follow from $(1-t^r)=\prod(1-t\zeta_r^k)$.

\bigskip
They yield: 
$$\prod_k \theta(x+{r \over r}-z)=
\prod_k q^{1 \over 8} 2 \sin \pi (x+{k \over r}-z)\prod_l (1-q^l)
\prod_l(1-q^le^{2 \pi i (x+{k \over r}-z)})
  (1-q^le^{2 \pi i -(x+{k \over r}-z)})
$$
$$ 
=q^{r \over 8} 2^r {1 \over {2^{r-1}}} \sin \pi r(x-z) (\prod_l
  (1-q^l))^r
\prod_l(1-q^{rl}e^{2 \pi i r(x-z)})
  (1-q^le^{2 \pi i r(x-z)})$$
$$=
{{\prod_l  (1-q^l))^{r}} \over {\prod_l(1-q^{lr})}} \theta(r \tau, r(x-z)).$$
If $\sigma_r$ is a cyclic permutation of $X^r$ then 
the fixed point set is the diagonal, 
the representation of 
$\sigma_r$ in the normal bundle is the quotient of the
regular representation by the trivial representation and each isotrivial 
component isomorphic to the tangent bundle of  $X$.
Therefore:
$$L(\sigma_r,X^r)=\prod_i 
\prod_{k=0}^{r-1} y_i {{\theta({y_i \over {2 \pi i}}+{k \over r}-z)}
\over {\theta({y_i \over {2 \pi i}}+{k \over r})}}[X]=
\prod_i y_i {{\theta (r\tau,ry_i-rz))} \over {\theta(r \tau, ry_i)}}[X]=$$
$${1 \over r^d} \prod_i ry_i {{\theta (r\tau,ry_i-rz))} \over {\theta(r
\tau, ry_i)}}[X]=Ell(r\tau,y^r)  $$
(the latter equality follows since replacing $y_i \rightarrow ry_i$
multiplies the degree $d$ component of the cohomology class evaluated
on $[X]$ by $r^d$). 

\par We can use arguments similar to those used in 
\cite{macdonald},\cite{zagier} and \cite{hirzhofer} to conclude the proof
of \ref{generating}.
We have  
$$\sum Ell_n(X^n/\Sigma_n)t^n=\sum_n [{1 \over {\vert \Sigma_n \vert }}
\sum_{g \in {\Sigma_n}} L(g,X^n)]t^n$$ 
where $L(g,X^n)$ is the holomorphic Lefschetz number of $g$ acting
on the bundle ${\cal E}ll(X)$
As usual, one can replace the summation with the summation over the 
set of conjugacy classes
since 
conjugate $g$ have 
isomorphic fixed point sets. The number of elements in a conjugacy class is 
${\vert G \vert} \over {\vert C(g) \vert }$, where $C(g)$ is the centralizer of $g$. 
Hence the latter sum can be replaced by 
 $\sum_n \sum_{\{g\} \in \Sigma_n}{{L(g,X^n)}
 \over {\vert C(g) \vert }}t^n$.
Each conjugacy class is specified by a partition of $n$ which has $a_i$
cycles of length $i$, so that $\sum ia_i=n$. Let $g_{a_1,...,a_r}$ be an 
element in such a conjugacy class. 
Change of the order of summation yields
$$\sum_{a_1,...,a_n,...} L(g_{a_1,..,a_n},X^n){1 \over {(a_1)! \cdot \cdot\cdot {a_n}! \cdot 
2^{a_2} 
\cdot \cdot \cdot n^{a_n}}} t^{a_1+2a_2+...+na_n}$$ 
\noindent 
since the number of elements in the conjugacy class corresponding to $(a_1,..,a_n)$ is 
${{n!}\over {a_1!...a_n! 2^{a_2} ... n^{a_n}}}$.
Next, the fixed point set of $g_{a_1,..,a_n}$ is 
$X^{a_1} \times ... \times X^{a_n}$. Using the multiplicativity of
Lefschetz numbers we obtain
$$\sum_{a_1,...a_n,...} 
{{\prod_i L(\sigma_i,X^i)^{a_i} t^{a_1+2a_2+...na_n}} \over {a_1!...a_n! 2^{a_2}...n^{a_n}}}=
\prod_i \sum_k {{L(\sigma_i,X^i)^kt^{ki}} \over {{k!}i^k}}.$$ 
The latter can be simplified to
$$\prod_k \exp({{L(\sigma_i,X)t^i} \over i})=
\exp(\sum_{i,m,l} {{c(m,l)q^{im}y^{il}t^i} \over i})=$$
$$\prod_{m,l} \exp(-c(m,l)log (1-tq^my^l))=\prod_{m,n} {1 \over 
{(1-tq^my^l)^{c(m,l)}}}.$$ 

\hfill{$\Box$

\noindent We shall mention  the following special cases of \ref{generating}:

\begin{coro} Let $\chi_y(X)=\sum_p \chi^py^p$. Then
$$\sum_n \chi_y(X^n/\Sigma_n)t^n=
\prod_p{1 \over {(1-t(-y)^p)^{(-1)^p\chi^p}}}.$$
\end{coro}
This follows from \ref{generating} since
$\chi_y(X)=Ell(X,q=0,-y)(-y)^{d \over 2}$ and in particular
if $l+{d \over 2}=p$ then $c(0,l)=(-1)^p\chi^p$.
Generating series for 
$\chi_y$ were also considered in \cite{zhou1}, \cite{zhou2}.

\begin{coro} (Macdonald, \cite{macdonald}) Let $e$ denote the topological 
Euler characteristic. Then
$$\sum_n e(X^n/\Sigma_n)t^n={1 \over (1-t)^{e(X)}}.$$
\label{macdnfla}
\end{coro}

\begin{coro} (D.Zagier, \cite{zagier}) Let $\sigma$ denote the 
signature of the intersection form in the middle dimension.
Then $$\sum_n \sigma(X^n/\Sigma_n)t^n=
{{(1+t)^{{\sigma(X) -e(X)} \over 2} 
\over  (1-t)^{{\sigma(X)+e(X)} \over 2}}}.$$ 
\label{zagierfla}
\end{coro}


\begin{thebibliography}{99}

\bibitem{AtiyahHirz}M.Atiyah, F.Hirzebruch, {\em  Spin-manifolds and group
actions}. 1970 Essays on Topology and Related Topics (M\'emoires d\'edi\'es 
\`a Georges de Rham)
pp. 18--28 Springer, New York

\bibitem{AtiyahSinger} M.Atiyah, I.Singer, 
{\em The index of elliptic operators. III.}
Ann. of Math. (2) 87 1968 546--604. 

\bibitem{atiyahsegal} M.Atiyah, G.Segal, {\em 
On equivariant Euler characteristics.} J.
Geom. Phys. 6 (1989), no. 4, 671--677.


\bibitem{bat.dual} V.~V. Batyrev,  {\em Dual polyhedra and mirror
symmetry for Calabi-Yau hypersurfaces in toric varieties},
J. Algebraic Geom. 3 (1994) 493-535.

\bibitem{Batyrev.Dais} V.~V. Batyrev, D.~I. Dais, 
{\em Strong McKay correspondence, string-theoretic Hodge numbers and
mirror symmetry}, Topology, 35 (1996), no. 4, 901-929.

\bibitem{Batyrev} V.~V. Batyrev, 
{\em Non-Archimedean integrals and stringy Euler numbers 
of log-terminal pairs}, J. Eur. Math. Soc. (JEMS), 1 (1999), 
no. 1, 5--33.


\bibitem{Bvertex} L. A. Borisov, {\em Vertex Algebras and Mirror
Symmetry}, \\ preprint math.AG/9809094.

\bibitem{BorGun} L. A. Borisov, P. E. Gunnells, {\em Toric
varieties and modular forms}, preprint math.NT/9908138.

\bibitem{borlibg1} L.A.Borisov, A.Libgober, 
{\em Elliptic Genera of Toric 
Varieties and Applications to Mirror Symmetry}, 
Invent. Math. 140 (2000), p.453-485.

\bibitem{borlibg2} L.A.Borisov, A.Libgober,{\em Singular elliptic genus.}
preprint math.AG/0007108.
 
\bibitem{Taubes} R.Bott, C.Taubes, {\em On the rigidity theorems of
Witten}, Journal of A.M.S., {\bf 2} (1989), 138-186.

\bibitem{bryl} J.L.Brylinski, {\em Representations of loop groups,
Dirac operators on
loop space, and modular forms}. Topology 29 (1990), no. 4, 461--480.

\bibitem{Chandra} K. Chandrasekharan, {\em Elliptic functions}, Fundamental
Principles of Mathematical Sciences, 281, Springer-Verlag, Berlin-New York,
1985.

\bibitem{CoxKatz} D. Cox, S. Katz, {\em  Mirror Symmetry and Algebraic
Geometry}, Mathematical Surveys and monographs, 68, AMS, 1999.

\bibitem{Danilov} V.~I. Danilov, {\em The Geometry of Toric Varieties},
Russian Math. Surveys, 33 (1978), 97-154.

\bibitem{dijkgraaf} R.Dijkgraaf, {\em Fields, Strings, Matrices and 
Symmetric Products}, hep-th/9912104.

\bibitem{DMVV} R. Dijkgraaf, G. Moore, E. Verlinde, H. Verlinde, 
{\em Elliptic genera of symmetric products and second quantized strings},
Comm. Math. Phys. 185 (1997), no. 1, 197--209.

\bibitem{EOTY} T. Eguchi, H. Ooguri, A. Taormina, S.-K. Yang, 
{\em Superconformal algebras and string compactification on manifolds with 
$SU(N)$ holonomy}, Nucl. Phys. B315 (1989), 193.

\bibitem{EZ} M. Eichler, D. Zagier, {\em The theory of Jacobi forms},
Progress in Mathematics, 55, Birkh\"auser Boston, Inc., Boston, Mass.,
1985.

\bibitem{gottell} G. Ellingsrud, L. G\"ottsche, M. Lehn, {\em
On the Cobordism Class of the Hilbert Scheme of a Surface}, 
preprint math.AG/990409.

\bibitem{Fulton} W. Fulton, {\em Introduction to toric varieties},
Princeton University Press, 1993.

\bibitem{gottsche} L.G\"ottsche, {\em Orbifold-Hodge numbers of Hilbert schemes.}
Parameter spaces (Warsaw, 1994), 83--87, Banach Center Publ., 36, 
Polish Acad. Sci., Warsaw, 1996.

\bibitem{Gritsenko} V. Gritsenko, {\em Elliptic genus of Calabi-Yau
manifolds and Jacobi and Siegel modular forms},
preprint math.AG/9906190, to appear in Algebra i Analiz (St. Petersburg
Math. Journal) 11:5 (1999).

\bibitem{Gritsenko2} V. Gritsenko, {\em Complex vector bundles and
    Jacobi forms}, preprint math.AG 9906191, (1999).

\bibitem{Hirz1} F. Hirzebruch, {\em Topological methods in Algebraic Geometry},
translated from German and Appendix One by R. L. E. Schwarzenberger. With a
preface  to the third English edition by the author and Schwarzenberger. 
Appendix Two by A. Borel. Reprint of the 1978 edition. Classics in Mathematics,
Springer-Verlag, Berlin, 1995. 

\bibitem{hirzhofer} F.Hirzebruch and T.H\"ofer {\em On the Euler number of an
orbifold}. Math. Ann. 286 (1990), no. 1-3, 255--260.


\bibitem{Hirz2} F. Hirzebruch, {\em Elliptic genera of level $N$ for 
complex manifolds}, Differential Geometric methods in Theoretical Physics
(Como 1987). K. Bleuer, M. Werner Editors, NATO Adv. Sci.Inst.Ser. 
C: Math.Phys. Sci; 250. Dordrecht, Kluwer Acad.Publ.,1988. 

\bibitem{hirz70} F.Hirzebruch, 
{\em Complex cobordisms and the elliptic genus}, 
Notes by Adrian Langer. Contemp. Math., 241, 
Algebraic geometry: Hirzebruch 70 (Warsaw, 1998), 9--20, 
Amer. Math. Soc., Providence, RI, 1999. 

\bibitem{modforms} F.Hirzebruch,T.Berger, R.Jung,
{\em Manifolds and modular forms}. With appendices by Nils-Peter Skoruppa 
and by Paul Baum. Aspects of 
Mathematics, E20. Friedr. Vieweg and Sohn, Braunschweig, 1992.


\bibitem{hohn.thesis} G.H\"ohn, {\em Komplexe elliptische Geschlechter und 
$S^1$-\"equivariante Kobordismustheorie}. Diplomarbeit, Bonn. August, 
1991.



\bibitem{Kac} V.~Kac, {\em Vertex algebras for beginners}, University 
Lecture Series, {10}, American Mathematical Society, Providence, RI,
1997.

\bibitem{KYY} T. Kawai, Y. Yamada, S.-K. Yang,
{\em Elliptic Genera and N=2 Superconformal Field Theory}, Nucl. Phys.
{ B414} (1994), 191-212.

\bibitem{kosn} C.Kosniowski, {\rm Generators of the $Z/p$ bordism ring. 
Serendipity.}
Math. Z. 149 (1976), no. 2, 121--130.


\bibitem{Krichever} I. Krichever, {\em Generalized elliptic genera and
  Baker-Akhiezer functions}, Math. Notes, 47 (1990), 132-142.

\bibitem{Landweber.Stong} P. S. Landweber, editor, {\em Elliptic curves
and modular forms in algebraic topology}, Lecture Notes in Math.,
{ 1326}, Springer, Berlin, 1988.

\bibitem{LW} A. Libgober, J. Wood, {\em Uniqueness of the complex structure on
K\"{a}hler manifolds of certain homotopy types}, J. Differential Geom.
{32} (1990),
no. 1, 139--154.

\bibitem{liuridg1} Kefeng Liu, {\em On modular invariance and rigidity 
theorems.} J. Differential
Geom. 41 (1995), no. 2, 343--396.

\bibitem{liumod} Kefeng Liu, {\em Modular invariance and 
characteristic numbers.} Comm. Math. Phys. 174 (1995), no. 1, 29--42

\bibitem{liuridg} Kefeng Liu, {\em On elliptic genera and theta-functions}.
Topology 35 (1996), no. 3, 617--640. 

\bibitem{Liu} Kefeng Liu, {\em Modular Forms and Topology},
  Contemp. Math. {\bf 193}, AMS, 1996. 

\bibitem{MSV} F. Malikov, V. Schechtman, A. Vaintrob,
{\em Chiral de Rham complex}, preprint alg-geom/9803041.

\bibitem{macdonald} I.Macdonald, {\em The Poincar\'e polynomial of a symmetric 
product.} Proc. Cambridge Philos. Soc. 58 1962 563--568.

\bibitem{Mumford} D.Mumford, {\em Tata lectures on theta. I},  with the
assistance of C. Musili, M. Nori, E. Previato and M. Stillman. Progress in
Mathematics, {\bf 28}, Birkh\"{a}user Boston, Inc., Boston, Mass., 1983.

\bibitem{Neumann}  C. D. D. Neumann, {\em The elliptic genus of Calabi-Yau $3$-
and $4$-folds, product formulae and generalized Kac-Moody algebras}, 
J. Geom. Phys., {\bf 29} (1999), no. 1-2, 5--12.

\bibitem{ochanine1} S.Ochanine, {\em
Sur les genres multiplicatifs d\'efinis par des int\'egrales
elliptiques.} Topology 26 (1987), no. 2, 143--151.

\bibitem{Oda} T. Oda, {\em Convex Bodies and Algebraic Geometry - An  
Introduction to the Theory of Toric Varieties}, Ergeb. Math. Grenzgeb.
(3), vol. 15, Springer-Verlag, Berlin, Heidelberg, New York, London,
Paris, Tokyo, 1988.

\bibitem{segal} G.Segal, Seminar Bourbaki,
{\em Elliptic cohomology (after Landweber-Stong, Ochanine,
Witten, and others).} S\'eminaire Bourbaki, Vol. 1987/88. 
Ast\'erisque No. 161-162, (1988), Exp.
No. 695, 4, 187--201 (1989). 

\bibitem{tomanoi} H.Tamanoi, {\em 
Elliptic genera and vertex operator super-algebras.}
Lecture Notes in Mathematics, 1704. Springer-Verlag, Berlin, 1999.

\bibitem{Totaro} B. Totaro, {\em Chern numbers of singular varieties and
elliptic homology}, preprint, University of Chicago, 
to appear in Annals of Mathematics.

\bibitem{witten93} E.Witten, {\em On Landau-Ginzburg description of 
$N=2$ minimal models,}  Int. J. Mod. Phys. A9 (1994) 4783.

\bibitem{witten88} E.Witten, {\em  The index of the Dirac operator in loop
space}, Lecture Notes in Mathematics, vol. 1326. 1988. pp. 161--181.

\bibitem{zagier} D.Zagier,
{\em Equivariant Pontrjagin classes and applications to
orbit spaces.}
Applications of the $G$-signature theorem to transformation groups, symmetric
products and number theory. Lecture Notes in Mathematics, 
Vol. 290. Springer-Verlag, Berlin-New York, 1972. viii+130 pp

\bibitem{zaslow} E.Zaslow, {\em Topological orbifold models and quantum 
cohomology rings.}
Comm. Math. Phys. 156 (1993), no. 2, 301--331

\bibitem{zhou1} Jian Zhou, {\em Delocalized equivariant coholomogy of 
symmetric products,} preprint math.DG/9910028.

\bibitem{zhou2} Jian Zhou, 
{\em Calculations of the Hirzebruch $\chi_y$ genera of symmetric products 
by the holomorphic Lefschetz formula,} preprint math.DG/9910029.

\end{thebibliography}
\end{document}